\documentclass[11pt]{amsart}

\usepackage{amssymb}
\usepackage{epsf}
\usepackage{a4wide}

\newtheorem{theorem}{Theorem}
\newtheorem{fact}{Fact}
\newtheorem{coro}{Corollary}
 
\newcommand{\KK}{\mathbb{K}\hspace{0.5pt}}
\newcommand{\kk}{\Bbbk}
\newcommand{\CC}{\mathbb{C}\hspace{0.5pt}}
\newcommand{\RR}{\mathbb{R}\hspace{0.5pt}}
\newcommand{\ZZ}{\mathbb{Z}}
\newcommand{\II}{\mathbb{I}}
\newcommand{\QQ}{\mathbb{Q}\hspace{0.5pt}}
\newcommand{\HH}{\mathbb{H}\hspace{0.5pt}}
\newcommand{\bs}[1]{\boldsymbol{#1}}
\newcommand{\OO}{\mathcal{O}}
\newcommand{\oo}{{\scriptstyle\mathcal{O}}}
\newcommand{\ooo}{{\scriptscriptstyle\mathcal{O}}}

\begin{document}

\title[A Note on Shelling]{A Note on Shelling}

\author{Michael Baake}
\address{Institut f\"{u}r Mathematik, Universit\"{a}t Greifswald,
Jahnstr.~15a, 17487 Greifswald, Germany}
\email{mbaake@uni-greifswald.de}
\urladdr{http://schubert.math-inf.uni-greifswald.de/}

\author{Uwe Grimm}
\address{Applied Mathematics Department, 
Faculty of Mathematics and Computing,
The Open University, Walton Hall, 
Milton Keynes MK7 6AA, UK}
\email{u.g.grimm@open.ac.uk}
\urladdr{http://mcs.open.ac.uk/ugg2/}

\begin{abstract} 
The radial distribution function is a characteristic geometric
quantity of a point set in Euclidean space that reflects itself in the
corresponding diffraction spectrum and related objects of physical
interest.  The underlying combinatorial and algebraic structure is
well understood for crystals, but less so for non-periodic
arrangements such as mathematical quasicrystals or model sets.  In
this note, we summarise several aspects of central versus averaged
shelling, illustrate the difference with explicit examples, and
discuss the obstacles that emerge with aperiodic order.
\end{abstract}
\maketitle


\section{Introduction}
 
One characteristic geometric feature of a discrete point set
$\varLambda\subset\RR^d$, which might be thought of as the set of
atomic positions of a solid, say, is the number of points of
$\varLambda$ on shells of radius $r$ around an arbitrary, but fixed
centre in $\RR^d$. Of particular interest are special centres, such as
points of $\varLambda$ itself, or other points that are fixed under
non-trivial symmetries of $\varLambda$. This leads to the so-called
{\em shelling structure}\/ of $\varLambda$. Here, we consider infinite
point sets only. In general, one obtains different answers for
different centres, and one is then also interested in the average over
all points of $\varLambda$ as centres, called the {\em averaged
shelling}.

The spherical shelling of lattices and crystallographic point sets
(i.e., periodic point sets whose periods span ambient space) is well
studied, and many results are known in terms of generating
functions. If $\varLambda$ is a lattice, the number of points on
spheres of radius $r$ centre $0$ ({\em central shelling}\/) is usually
encapsulated in terms of the lattice theta function \cite[Ch.\
2.2.3]{CS}
\begin{equation} \label{theta1}
   \varTheta^{}_{\varLambda} (z) \; = \; 
   \sum_{x\in\varLambda} q^{x\cdot x}
         \; = \; \sum_{k} c(k)\, q^k
\end{equation}
where $q=e^{\pi i z}$ and $c(k)$ is the number of lattice points of
Euclidean square norm ($=$ square length) $k$. A closed expression for
the latter can be given in many cases, see \cite[Ch.\ 4]{CS} for
details on root and weight lattices and various related packings, and
\cite{BGJR} for an explicit example. There are many related lattice
point problems, see \cite{Hux} and references therein for recent
developments.

One special feature of a lattice is that the shelling generating
function is independent of the lattice point which is chosen as the
centre -- and, consequently, {\em central}\/ and {\em averaged}\/
shelling give the same result. Similarly, for uniformly discrete point
sets that are crystallographic, the average is only over finitely many
points in a fundamental domain and can often be calculated
explicitly. For general uniformly discrete point sets, however, the
situation is more complicated in that no two centres might give the
same shelling function, or that the average may not be well
defined. But there is one important class of point sets, the so-called
{\em model sets}\/ (also called cut-and-project sets, see
\cite{Moody,MB,Martin1,S,LM,BS,BSW} and references therein), which
provide a high degree of order and coherence so that an extension of
the shelling problem to these cases is possible, and has indeed been
pursued. The original motivation for the investigation of model sets
came from applications in physics. Meanwhile, due to interesting
connections with several branches of mathematics, they are also
studied in their own right, see \cite{Mbook,Pat,BMbook} and
\cite{BGlit} for details and further references.  Below, we shall
summarise the key properties of model sets needed for this
article.\smallskip

One of the earliest attempts to the shelling of model sets, to our
knowledge, is that of Sadoc and Mosseri \cite{SM} who investigated the
4D Elser-Sloane quasicrystal \cite{ES} and then conjectured a formula
for the central shelling of a close relative of it which was obtained
by replacing the highly symmetric 4D polytype used in \cite{ES} by a
4D ball as a window. The conjecture was put right and proved in
\cite{MW} by means of algebraic number theory revolving around the
arithmetic of the icosian ring $\II$, a maximal order in the
quaternion algebra $\HH(\QQ(\mbox{\small $\sqrt{5}\,$}))$. Recently,
the central shelling was extended to the much more involved 3D case of
icosahedral symmetry \cite{Al}.  Also, some results exist on planar
cases, e.g., for special eightfold and twelvefold symmetric cases with
circular windows \cite{MS,MS2}.

The common aspect of all these extensions to model sets (or
mathematical quasicrystals) is that only the central shelling of a
highly symmetric representative has been considered, with a ball as
window in internal space. This is a rather special situation which
appears slightly artificial in view of the fact that the most relevant
and best studied model sets usually have polytopes rather than balls
as window, or, more generally, even compact sets with fractal
boundary, cf.\ \cite{MB}. Also mathematically, the classical examples
such as the rhombic Penrose or the Ammann-Beenker tiling are very
attractive due to their rather intricate and unexpected topological
nature \cite{AP,FHK}. 

A more natural approach to model sets seems to be the {\em averaged}\/
shelling, and it is the aim of this article to start to develop this
idea. As we shall see, the topological structure will be manifest in
the examples discussed below.  On the other hand, the {\em central}\/
shelling does have a universal meaning, too, if one considers it first
for modules rather than for model sets. The window condition can then
be imposed afterwards, see \cite{BGJR,MU1,MU2} for some examples.
This approach is implicit in \cite{MW}, but does not seem to have
attracted much notice. It is important, though, because it leads to a
separation of universal and non-universal aspects.

\section{Central shelling}

Many of the well studied planar tilings with non-crystallographic
symmetries share the property that their vertices (or other typical
point set representatives) form a discrete subset of rings of
cyclotomic integers. This gives a nice and powerful link to results of
algebraic number theory, which has, in fact, been used to construct
model sets \cite{P}, and which also appeared before in a different
context \cite{Mermin}. Let us thus first explain the situation of
central shelling for these underlying dense point sets.\smallskip

Let $\xi_{n}$ be a primitive $n$-th root of unity (with $n\geq 3$),
e.g., $\xi_{n}=e^{2\pi i/n}$, and $\QQ(\xi_{n})$ the corresponding
cyclotomic field. Then, $\QQ(\xi_{n}+\overline{\xi}_{n})$ is its maximal
real subfield. From now on, we will use the following notation
\begin{equation}
\KK=\QQ(\xi_{n})\,,\quad
\kk=\QQ(\xi_{n}+\overline{\xi}_{n})\,,\quad
\OO = \ZZ[\xi_{n}]\,,\quad
\oo = \ZZ[\xi_{n}+\overline{\xi}_{n}]\,,
\label{nota}
\end{equation} 
where $\OO$ is the ring of cyclotomic integers, which is the maximal
order of $\KK$, and $\oo$ is the ring of algebraic integers of $\kk$,
see \cite{Wash}.

Note that $\OO$ is a $\ZZ$-module of rank $\varphi(n)$, where
$\varphi$ denotes Euler's totient function. The set $\OO$, seen as a
(generally dense) point set in $\RR^2$, has $N$-fold rotational
symmetry, where
\begin{equation}
N\; =\; N(n)\;=\;
    \begin{cases}
         n  & \text{if $n$ is even,} \\
         2n & \text{if $n$ is odd.}
    \end{cases}
\label{symm}
\end{equation}
This also means that $\OO$ has precisely $N$ units on the unit circle,
which are actually all roots of unity of $\KK$. Also, $\KK$ is a totally
complex field extension of $\kk$ of degree 2. It is known that, in this
cyclotomic situation, the unique prime factorisation property of $\OO$
(i.e., class number one) implies that of $\oo$, and this happens in
precisely 29 cases, compare \cite[Thm.~11.1]{Wash}, namely for
\begin{equation}
\begin{split}
n \in \{ &3,4,5,7,8,9,11,12,13,15,16,17,19,20,21,\\
         &24,25,27,28,32,33,35,36,40,44,45,48,60,84\}\,,
\end{split}
\label{list}
\end{equation}
where $n\not\equiv 2\bmod 4$ to avoid double counting.  Note that
$n=1$ ($N=2$) is excluded here because it corresponds to $\KK=\QQ$ with
$\OO=\ZZ$, which is only one-dimensional.

Now, let $p$ be a prime of $\oo$. Then, in going from $\oo$ to $\OO$,
precisely one of the following cases applies, see \cite[Ch.~I]{Neu}
and \cite[Ch.~4]{Wash}:
\begin{itemize}
\item[(1)] $p$ {\em ramifies}, i.e., $p=P\overline{P}$ with $P$ a prime and
$\overline{P}/P$ a root of unity in $\OO$.
\item[(2)] $p$ is {\em inert}, i.e., $p$ is also prime in $\OO$.
\item[(3)] $p$ is a {\em splitting prime} of $\OO/\oo$, i.e.,
           $p=P\overline{P}$ with $\overline{P}/P$ not a unit in $\OO$.
\end{itemize}
Up to units, all primes of $\OO$ appear this way.

Prime factorisation in $\oo$ versus $\OO$ can now be employed to find
the combinatorial structure of the shells. We encode this into the
{\em central shelling function}\/ $c(r^2)$ which counts the number of
points on shells (circles) of radius $r$. By convention, $c(0)=1$.

\begin{theorem} \label{thm1}
  Let\/ $\OO=\ZZ[\xi_{n}]$ be any of the\/ $29$ planar\/ $\ZZ$-modules
  that consist of the integers of a cyclotomic field with class number
  one.  Then, for\/ $r^{2}>0$, the function\/ $c(r^2)$ vanishes
  unless\/ $r^{2}\in\oo$ and all inert prime factors of\/ $r^2$ occur
  with even powers only. In this case, 
\begin{equation} \label{shell-for}
    c(r^2) \; = \; N \cdot \!\!\prod_{\stackrel{\scriptstyle p\mid r^2}
            {\scriptstyle p \; {\rm splits}}} \!\!\big(t(p)+1\big)\, ,
\end{equation} 
  where\/ $p$ runs through a representative set of the primes of\/
  $\oo$. Here, $t(p)$ is the maximal power\/ $t$ such that\/ $p^t$
  divides\/ $r^2$.  The prefactor, $N=N(n)$ of Eq.~$(\ref{symm})$,
  reflects the point symmetry of the module. Furthermore, $r^2$ is
  then a totally positive number in\/ $\oo$, i.e., all its algebraic
  conjugates are positive as well. 
\end{theorem}

\begin{proof}

Since $c(0)=1$ by convention, consider $r^{2}>0$. If there exists a
number $x\in\OO$ on the shell of radius $r$ around $0$, we must have
$r^2=x\overline{x}$, hence $r^2\in\oo$. In this case, any inert prime
factor $p$ of $r^2$ (in $\oo$) necessarily divides both $x$ and
$\overline{x}$ (in $\OO$). Consequently, the maximal power $t=t(p)$
such that $p^t$ divides $r^2$ must be even.

Conversely, assume $r^{2}>0$ and $t(p)$ even for all inert primes of
$\oo$. If a ramified or a splitting prime $p=P\overline{P}$ divides $r^2$,
we know that equal powers of $P$ and $\overline{P}$ occur in the prime
factorisation of $r^2$ in $\OO$. Consequently, we can group the prime
factors of $r^2$ in $\OO$ into two complex conjugate numbers, i.e., we
have at least one solution of the equation $r^2=x\overline{x}$ with
$x\in\OO$, so $c(r^2)>0$.

Consider a non-empty shell with $r^{2}>0$, i.e., $r^2=x\overline{x}$
for some $0\ne x\in\OO$.  Consider the prime factorisation $r^2 =
e\cdot p_1^{t^{}_1}\cdot\ldots\cdot p_s^{t^{}_s}$ in $\oo$, with $e$ a
unit. If $p_i$ is {\em not}\/ a splitting prime, the distribution of
the corresponding primes in $\OO$ to $x$ and $\overline{x}$ is unique,
up to units of $\OO$.

If, however, $p_j = P_j \overline{P}_j$ {\em is}\/ a splitting prime, we
have to distribute $(P_j \overline{P}_j)^{t^{}_j}$ over $x$ and
$\overline{x}$. In view of $\overline{P}_j$ being the complex conjugate of
$P_j$, but not an algebraic conjugate, we have the options of
$(P_j)^s_{}(\overline{P}_j)^{t_j - s}_{}$ as factor of $x$ and
$(\overline{P}_j)^s_{}(P_j)^{t_j - s}_{}$ as factor of $\overline{x}$, for any
$0\le s\le t^{}_j$. This amounts to $t^{}_j +1$ different
possibilities, which gives the corresponding factor in
(\ref{shell-for}).

As mentioned above, there are $N$ units of $\OO$ on the unit
circle. This means that, as soon as $r^2>0$, points on the shells come
in sets of $N$, which gives the prefactor in
(\ref{shell-for}). Together with the previous arguments, this explains
the multiplicative structure of $c/N$.

Finally, assume $r^2=x\overline{x}$ for some $0\ne x\in\OO$ and let
$\sigma$ be any Galois automorphism of $\KK$ over $\QQ$. Then we have
\[
0\; <\; \sigma(x)\overline{\sigma(x)}\; =\; \sigma(x)\sigma(\overline{x})
\;=\; \sigma(x\overline{x}) \;=\; \sigma(r^2)
\]
so also all algebraic conjugates of $r^2$ are positive. This shows
that $r^2$ is totally positive.
\end{proof}

\smallskip
\noindent
{\sc Remark}: It is clear that Theorem \ref{thm1} can be generalised
to the situation that $\KK$ is a totally complex field extension of a
totally real field $\kk$ whenever $\KK$ has class number one, with sets of
integers $\OO$ and $\oo$ as above, compare \cite[Thm.~4.10]{Wash}. In
this case, the prefactor in Eq.~(\ref{shell-for}) has to be replaced
by the number of elements in the unit group of $\OO$ that lie on the
unit circle.\smallskip

Let us consider the cyclotomic case in more detail. If $V(r^2)=\{x\in
\OO\mid x\overline{x}=r^2\}$, and $\sigma$ is any Galois automorphism of
$\KK/\QQ$, then $\sigma(\OO)=\OO$ and $V(r^2)$ is mapped bijectively to
$V(\sigma(r^2))$. This means that $c(r^{2})=c(\sigma(r^2))$.

Moreover, consider the situation that two totally positive numbers of
$\oo$, $r^2$ and $R^2$, are related by $R^2=er^2$, with $e$ a unit in
$\oo$. Clearly, $e$ is then also totally positive. If $e$ is of the
form $e=u\overline{u}$, with $u$ a unit in $\OO$, the mapping $x\mapsto ux$
gives a bijection between $V(r^2)$ and $V(R^2)$, hence
$c(r^2)=c(R^2)$.  If all totally positive units of $\oo$ are of this
form, which includes the case that $e$ is the square of a unit in
$\OO$, we may conclude that the central shelling function $c$ only
depends on the principal ideal of $\oo$ generated by $r^2$.

In general cyclotomic fields, this factorisation property of totally
positive units need not be satisfied (e.g., it fails for
$n=29$). However, it is true for all class number one cases. More
precisely, if $n$ is a power of $2$, all totally positive units of
$\oo$ are squares of units of $\OO$, which is known as Weber's
theorem, compare \cite[Cor.~1 and Rem.~2]{Gar}.  The same statement
holds if $n$ is an odd prime below 100, except for $n=29$, see
\cite[Ex.~2]{Gar}. We checked explicitly, using the KANT program
package \cite{DFKPRW,kant}, that this remains true for all $n$ from
our list (\ref{list}) that are prime powers. All remaining cases of
(\ref{list}) are composite integers. Here, not all totally positive
units of $\oo$ are squares in $\OO$, but they are of the form
$e=u\overline{u}$, with $u$ a unit in $\OO$. This was again checked
using KANT. The difference to the other cases comes from the
additional unit $u=1-\xi_{n}$, compare \cite[Cor.~4.13]{Wash}.

We may conclude as follows. 
\begin{fact}\label{fact1}
Let\/ $\OO=\ZZ[\xi_{n}]$ be any of the\/ $29$ planar\/ $\ZZ$-modules
that consist of the integers of a cyclotomic field with class number
one, and\/ $\oo=\ZZ[\xi_n+\overline{\xi}_n]$. Then, the central shelling
function\/ $c$ for\/ $\OO$ depends on $r^2\in\oo$ only via the
principal ideal\/ $r^2\oo$ generated by it.\hfill\qed
\end{fact}

This allows us to reformulate the result of Theorem~\ref{thm1} by
means of ideals and characters of the field extension $\KK/\kk$. By a
character $\chi\not\equiv 0$, we here mean a totally multiplicative
real function of the ideals of $\oo$, i.e., $\chi(\mathfrak{ab})=
\chi(\mathfrak{a})\chi(\mathfrak{b})$ for all ideals $\mathfrak{a}$
and $\mathfrak{b}$ of $\oo$, see \cite[Ch.~VII.6]{Neu} for background
material.  In particular, $\chi(\oo)=1$. It suffices to specify the
values of $\chi$ for all prime ideals $\mathfrak{p}$ of $\oo$. We
define
\begin{equation}
\chi(\mathfrak{p})\;=\;
\begin{cases}
\hphantom{-}0 & \text{if $\mathfrak{p}$ ramifies}\\
-1 & \text{if $\mathfrak{p}$ is inert}\\
\hphantom{-}1 & \text{if $\mathfrak{p}$ splits}
\end{cases}
\label{chardef}
\end{equation}
where the property of the prime ideal $\mathfrak{p}$ refers to the
behaviour under the field extension from $\kk$ to $\KK$.  This leads to
the following result.
\begin{coro}\label{cor1}
Under the assumptions of Theorem~$\ref{thm1}$, the central shelling
function\/ $c$ is proportional to the summatory function of the
character\/ $\chi$ of Eq.~$(\ref{chardef})$, i.e.,
\begin{equation}
c(r^2\oo)\;=\;
N\cdot\sum_{\mathfrak{a}\mid (r^{2}\ooo)}\chi(\mathfrak{a}),
\label{ideal}
\end{equation}
with\/ $N$ given by Eq.~$(\ref{symm})$.
\end{coro}

\begin{proof}
Due to unique prime factorisation in $\oo$ and the multiplicative
structure of $c/N$ according to Eq.~(\ref{shell-for}), it is
sufficient to verify the claim for prime powers, i.e., for
$r^2\oo=\mathfrak{p}^{\ell}$. Clearly, if $\mathfrak{p}$ ramifies, the
sum in Eq.~(\ref{ideal}) gives $c(\mathfrak{p}^{\ell})=1$ for all
$\ell\ge 0$. The alternating sign of $\chi(\mathfrak{p}^\ell)$ for
inert $\mathfrak{p}$ implies $c(\mathfrak{p}^\ell)=0$ for odd $\ell$
and $c(\mathfrak{p}^{\ell})=1$ otherwise.  If $\mathfrak{p}$ splits,
the right hand side of Eq.~(\ref{ideal}) adds up to $\ell+1$. Invoking
Fact~\ref{fact1} and a comparison with Eq.~(\ref{shell-for}) completes
the proof.
\end{proof}
\smallskip

The explicit use of Theorem~\ref{thm1} and Corollary~\ref{cor1}
requires the knowledge of the splitting structure of the primes.
Examples can be found in \cite{PBR,MU1,MU2}, see also \cite{MS,MS2}.
If one is interested in the central shelling of a {\em model set\/}
rather than that of the underlying (dense) module, one has to take the
window into account as a second step. A model set
$\varLambda(\varOmega)$ in ``physical space'' $\RR^d$ is defined
within the following cut-and-project scheme \cite{Moody,MB}
\begin{equation}
\renewcommand{\arraystretch}{1.5}
\begin{array}{ccccc}
\RR^{d} & \xleftarrow[\qquad]{\pi} & \RR^{d}\times H & 
\xrightarrow[\qquad]{\pi_{H}^{}} & H \\
\cup && \cup && \cup \makebox[0pt][l]{\small\ \ dense}\\ 
L & \xleftarrow[\qquad]{1-1} & \varGamma & \xrightarrow[\qquad]{} & L^*\\
\end{array}
\label{candp}
\end{equation}
where the ``internal space'' $H$ is a locally compact Abelian group,
and $\varGamma\subset\RR^{d}\times H$ is a lattice, i.e., a co-compact
discrete subgroup.  The projection $L^*=\pi_{H}^{}(\varGamma)$ is
assumed to be dense in internal space, and the projection into
physical space has to be one-to-one on $\varGamma$. Consequently, the
mapping ${\,}^*$: $L\longrightarrow L^{*}\subset H$, with
${\,}^*=\pi_{H}^{}\circ\big(\pi|_{\varGamma}^{}\big)^{-1}$, is well
defined.  It is called the $*$-map of the cut-and-project formalism,
compare \cite{Moody}. Note that the $*$-map need not be injective,
i.e., its kernel can be a nontrivial subgroup of $L$.

A model set $\varLambda(\varOmega)$ is now defined as
\begin{equation}
   \varLambda(\varOmega)\;=\;
    \left\{x\in L\mid x^{*}\in\varOmega\right\}
   \;=\;
   \left\{\pi(y)\mid y\in\varGamma,\,
   \pi_{H}^{}(y)\in\varOmega \right\}\;\subset\; \RR^{d},
   \label{eq:ms}
\end{equation}
where the window $\varOmega\subset H$ is a relatively compact set with
non-empty interior. Usually, one either takes an open set or a compact
set that is the closure of its interior. Note that the $*$-map is
well defined on $\varLambda(\varOmega)$, with
$\big(\varLambda(\varOmega)\big)^* \subset\varOmega$. More
generally, also sets of the form $t+\varLambda(\varOmega)$ with
$t\in \RR^d$ are called model sets. If $t\in L$, one has
$t+\varLambda(\varOmega) = \varLambda(t^* + \varOmega)$ and is
back to the case of Eq.~(\ref{eq:ms}), which is sufficient for
our discussion. For the above example of a
cyclotomic field $\KK=\QQ(\xi_n)$, we need $d=\varphi(n)$ to construct
model sets with $n$-fold symmetry, compare \cite[App.~A]{BJS}.

In order to compute the central shelling for a model set
$\varLambda(\varOmega)$, one first determines all points of the module
$L=\pi(\varGamma)$ on the shell of a given radius $r$. Then, the
window $\varOmega$ decides, according to the filtering process of
Eq.~(\ref{eq:ms}), which of these points actually appear in the model
set, and the shelling formula is modified accordingly. As long as we
are dealing with a one-component model set (i.e., as long as all
points are in one translation class), the formula of Theorem
\ref{thm1} thus gives an upper bound on the shelling number in the
model set. As mentioned above, the central shelling of a few model
sets with spherical windows \cite{SM,MW,MS,MS2,Al} has been considered
in detail.

\section{Averaged shelling}

A moment's reflection reveals that the averaged shelling is
considerably more involved. In order to determine the averages, one
would need to know all possible local configurations up to a given
diameter together with their frequencies, provided the latter are well
defined. In general, this is not the case, as cluster or patch
frequencies in general Delone sets need not exist. However, regular
model sets are particularly nice in this respect because all patch
frequencies exist uniformly \cite{Martin1}, which is equivalent to
unique ergodicity of the corresponding dynamical system
\cite{Solo,Martin2} (under the translation action of
$\RR^d$). Moreover, due to existence of the cut-and-project scheme
(\ref{candp}) and Weyl's theorem, compare \cite{Moody2001}, it is
possible to transfer the averaging part of the combinatorial problem
to one of analysis.

Let us also point out that Eq.~(\ref{theta1}) for a model set does not
make much sense as it would depend on the representative chosen,
rather than being a quantity attached to an entire local
indistinguishability (LI) class, compare \cite{S,MB}. If $\varLambda$
is a lattice, $\varLambda-\varLambda=\varLambda$, and we could equally
well sum over the difference set in (\ref{theta1}).  Using this for
model sets would give $\sum_{x\in \varLambda-\varLambda} q^{x\cdot x}$
which is constant on the LI class. However, this still does not
reflect the {\em statistical}\/ aspects of the (local) shells, because
each $x\in\varLambda-\varLambda$ is counted with weight one. Let us
thus introduce the {\em averaged shelling function}\/ $a(r^2)$ as the
number of points on a shell of radius $r$, averaged over all points of
$\varLambda$ as possible centres of the shells.

Now, let $\varLambda = \varLambda(\varOmega)$ be a regular, generic
model set, in the terminology of \cite{Moody}, with window
$\varOmega$, i.e., $\varOmega$ is a relatively compact set in $H$ with
non-empty interior, boundary of measure $0$, and
$\partial\varOmega\cap\pi_{H}^{}(\varGamma)=\varnothing$. For
simplicity, we also assume that $H=\RR^{m}$, though a generalisation
of what we say below to the case of general locally compact Abelian
groups is possible. In analogy to Eq.~(\ref{theta1}), a generalised
theta series could be defined ad hoc as
\begin{equation} \label{theta2}
   \varTheta^{}_{\varLambda}(z) \; := \; 
         \sum_{r\in\mathcal{R}} a(r^2)\, q^{r^2}
\end{equation}
where $q=e^{\pi iz}$ and $\mathcal{R} = \{ r\in\RR_{\geq 0}\mid |y|=r
\mbox{ for some } y\in \varLambda-\varLambda\}$ is the set of possible
radii as obtained from the set of difference vectors between points of
$\varLambda$. The coefficient $a(r^2)$ is now meant as the averaged
quantity defined above, which we will now calculate.

Let $\nu(y)$ denote the relative frequency of the difference $y$
between two points of the model set (hence
$y\in\varLambda-\varLambda$). Up to the overall density of the model
set, $\nu(y)$ is an autocorrelation coefficient of the point set
$\varLambda$.  This quantity exists uniformly for all $y$ as a
consequence of the model set structure
\cite{Hof,Martin1,Moody2001}. But then, we obviously obtain
\begin{equation} \label{ashell1}
    a(r^2) \; = \; 
      \sum_{\stackrel{\scriptstyle y\in\varLambda-\varLambda}
{\scriptstyle |y|=r}} 
            \, \nu(y)\, .
\end{equation}
On the other hand, if $\varLambda_s=\{x\in\varLambda\mid |x| < s\}$,
one has
\begin{eqnarray} \label{weyl}
  \nu(y) & = & \lim_{s\to\infty} \frac{1}{|\varLambda_s|}
     \sum_{\stackrel{\scriptstyle x\in\varLambda_s}
{\scriptstyle x+y\in\varLambda}} 1
     \;\, = \;\, \lim_{s\to\infty} \frac{1}{|(\varLambda_s)^*|}
     \sum_{\stackrel{\scriptstyle x^*\in(\varLambda_s)^*}
                   {\scriptstyle (x+y)^*\in\varOmega}} 1 \nonumber \\
  & = & \frac{1}{{\rm vol}(\varOmega)}
     \int_{\RR^m} \bs{1}^{}_{\varOmega}(z)\, 
     \bs{1}^{}_{\varOmega}(z + y^*) \, {\rm d}z
\end{eqnarray}
where $\bs{1}^{}_{\varOmega}$ is the characteristic function of the
window. Note that, as the $*$-map need not be injective, the second 
equality may only hold in the limit $s\to\infty$ (this step is implicit
in the proof of \cite[Thm.~1]{Martin1}). We add it here because
it shows how the counting is transfered to internal space, in particular 
in the cases where the $*$-map {\em is\/} one-to-one, which is the
situation we will meet in the examples.

The last step in (\ref{weyl}) is now a direct application of Weyl's 
theorem on uniform distribution.  This is justified here because
$(\varLambda_s)^*$, for increasing $s$, gives a sequence of points in
$\varOmega$ that are uniformly distributed, see \cite{Hof,Moody2001}
and \cite[Thm.~1]{Martin1}, and because $\partial\varOmega$ has
measure $0$ by assumption. In this situation, the averaged quantities
are the same for generic and singular members of the LI class
\cite{Martin1,MB}. Moreover, it also does not change if
$\,\partial\varOmega\cap L^*\ne\varnothing$, so that the corresponding
assumption can be dropped.  Consequently, the averaged shelling
function is constant on LI classes of regular model sets. We combine
Eqs.~(\ref{ashell1}) and (\ref{weyl}) to obtain
\begin{theorem} \label{thm2}
  Let\/ $\varLambda$ be a regular model set in the sense of
  Moody~\cite{Moody}, obtained from a cut-and-project scheme\/
  $(\ref{candp})$ with internal space\/ $H=\RR^m$ and window\/
  $\varOmega$.  Then, the averaged shelling function\/ $a(r^2)$ exists,
  and is given by
\begin{equation}
\label{ashell2}
    a(r^2) \; = \; \frac{1}{{\rm vol}(\varOmega)}
      \sum_{\stackrel{\scriptstyle y\in\varLambda-\varLambda}
{\scriptstyle |y|=r}} 
   {\rm vol}\big(\varOmega\cap (\varOmega-y^*)\big).
\end{equation}
In particular, $a(r^2)$ vanishes if there is no\/
$y\in\varLambda-\varLambda$ with\/ $y\cdot y=r^2$.\qed
\end{theorem}

\noindent
{\sc Remark}: This result allows the calculation of the shelling
function, for any possible radius $r$, by evaluating {\em finitely}\/
many volumes in internal space. This is so because a model set
$\varLambda$ has the additional property that also its difference set,
$\varLambda - \varLambda$, is uniformly discrete, so that there are
only finitely many different solutions of $|y|=r$ with
$y\in\varLambda-\varLambda$.

\section{Examples}

Let us first consider a well-known model set in one dimension, the
Fibonacci chain, which can be described as
\begin{equation}
   \varLambda_{\text{F}} \; = \; 
\big\{ x\in\ZZ[\tau] \mid x^*\in [-1,\tau -1]\,\big\}
\; = \; \varLambda\big([-1,\tau-1]\big),
\label{fibo}
\end{equation}
where $\ZZ[\tau] = \{ m+n\tau\mid m,n\in\ZZ\}$ is the ring of integers
in the quadratic field $\QQ(\tau)$ and $\tau=(1+\sqrt{5}\,)/2$ is the
golden ratio.  The $*$-map in this setting is algebraic conjugation in
$\QQ(\tau)$, defined by $\sqrt{5}\mapsto -\sqrt{5}$. The 2D lattice
behind this formulation is $\varGamma=\{(x,x^*)\mid x\in\ZZ[\tau]\}$.
A short calculation results in
$\varLambda_{\text{F}}-\varLambda_{\text{F}}=
\varLambda\big([-\tau,\tau]\big)$, and
\begin{equation}\label{fibshell}
   \nu(y) \; = \; \nu(-y) \; = \; f^{}_{\text{F}}(y^*)
    \; = \; 
    \begin{cases}
         0              & \text{if $|y^*| >  \tau$} \\
         1 - |y^*|/\tau & \text{if $|y^*|\leq\tau$}
    \end{cases}
\end{equation}
so that the averaged shelling function for the Fibonacci chain (and
thus also for its entire LI class) is $a(0)=1$ and $a(r^2)=2
f^{}_{\text{F}}(y^*)$ for any non-zero distance $r$ that is the
absolute value of a number $y\in\varLambda_{\text{F}}
-\varLambda_{\text{F}} \subset\ZZ[\tau]$. Also, all shelling numbers
$a(r^2)$ are elements of $\ZZ[\tau]$, as can easily be seen from
formula (\ref{fibshell}). This has a topological interpretation, as we
will briefly explain below for a more significant example.

In internal space, the function $f^{}_{\text{F}}$ has a piecewise
linear continuation, but the function $a(r^2)$ looks rather erratic,
compare \cite{BGJR} for a similar example.  This is a consequence of
the properties of the $*$-map, being algebraic conjugation in this
case. As a mapping, it is totally discontinuous on $L^*$ (and also
on its rational span) when the latter is given the induced
topology of the ambient space $H$. In a different topology,
however, this map becomes uniformly continuous, and it is this
alternative setting, compare \cite{BM2001}, which explains the
appearance of the internal space from intrinsic data of a model set
$\varLambda$.

\begin{figure}[b]
\centerline{\epsfxsize=0.6\textwidth\epsfbox{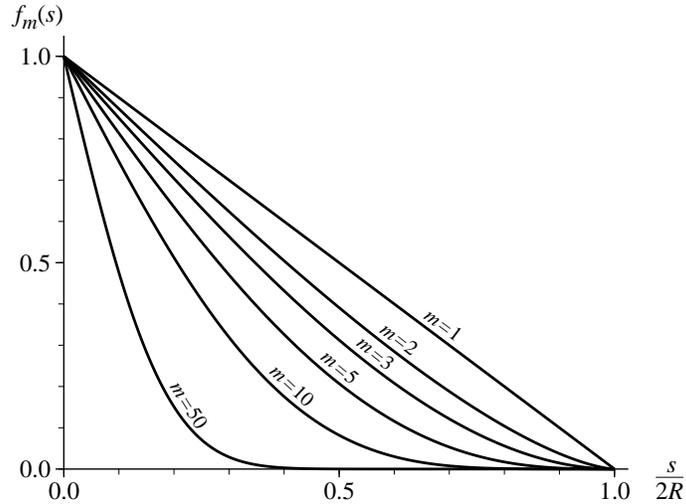}}
\caption{Radial component of the frequency functions $f^{}_{m}(s)$ of
Eqs.~(\ref{eq:f2}), (\ref{eq:feven}) and (\ref{eq:fodd}) for dimensions
$m=1,2,3,5,10,50$ of internal space.\label{fig:rad}}
\end{figure}

\smallskip
As another example, let us once more look at the circular shelling in
the plane, i.e., at a 2D model set with an open disk (radius $R$,
centre $0$) as window in 2D internal space. So, $\varOmega=B_R^{}(0)$
and, consequently, $\varOmega-\varOmega=B_{2R}^{}(0)$.  We have
$\nu(y)=f^{}_{2}(y^*)$, where, due to rotational symmetry of the
window, the function $f^{}_2$ only depends on $s=|y^*|$. Explicitly,
it is given by
\begin{equation}
  f^{}_{2}(s) \; = \; 
\frac{{\rm vol}\big(B_R^{}(0) \cap B_R^{}(s)\big)}
{{\rm vol}\big(B_R^{}(0)\big)} \; = \;
       \frac{2}{\pi}\, \arccos\left(\frac{s}{2R}\right) \; - \;
       \frac{s}{\pi R} \sqrt{1 - \left(\frac{s}{2R}\right)^{2}}
\label{eq:f2}
\end{equation}
for $0\le s<2R$ and $f^{}_2(s)=0$ otherwise. Fig.~\ref{fig:rad}
contains a graph of $f^{}_{2}(s)$. This function, often called the
{\em covariogram}\/ of the disk, is a radially symmetric positive
definite function known as Euclid's hat, see \cite[p.~100]{Gneit}.

To calculate $a(r^2)$, one has to sum finitely many terms of this
kind, according to Eq.~(\ref{ashell1}).  This situation of a 2D
internal space shows up for planar model sets with $n\in\{5,8,12\}$,
because these are the cases with $\varphi(n)=4$. Here, one simply
obtains
\begin{equation}
a(r^2) \; = \; c(r^2)\, f^{}_{2}(s)
\label{a-c}
\end{equation}
where $c(r^2)$ is the central shelling function of
Eq.~(\ref{shell-for}) and $s=|y^{*}|$ for any $y$ on the shell of
radius $r$. This is so because the window is a disk and the $*$-map
sends all cyclotomic integers on a circle to a single circle in
internal space. Consequently, the central shelling provides an upper
bound for the average shelling in this case.\smallskip

Part of this result can be extended to arbitrary dimension. For two
intersecting $m$-dimensional balls of radius $R$, the overlap consists
of two congruent ball segments. The corresponding volume can be
calculated by integrating slices (which are balls of dimension $m-1$).
Dividing by the volume of the $m$-ball, the covariogram becomes
\begin{equation}
  f^{}_{m}(s) \; = \;
  \frac{2\,\Gamma(\frac{m}{2}+1)}{\sqrt{\pi}\,\Gamma(\frac{m+1}{2})}
  \int\limits_{0}^{\arccos(\frac{s}{2R})} \sin^{m}(\alpha)\, 
  {\rm d}\alpha.
\label{eq:fgen}
\end{equation}
The integral can be expanded in terms of Chebyshev polynomials.  For
even $m=2\ell$, this yields
\begin{equation}
\label{eq:feven}
\begin{split}
f^{}_{2\ell}(s) \;=\;{}&
    \frac{\Gamma(\ell+1)}{2^{2\ell-1}\sqrt{\pi}\,\Gamma(\ell+\tfrac{1}{2})}
    \Bigg[\binom{2\ell}{\ell}\,\arccos\left(\frac{s}{2R}\right)\\[1mm]
    & +\; \sqrt{1-\left(\frac{s}{2R}\right)^{2}}\,
    \sum\limits_{k=1}^{\ell}\frac{(-1)^k}{k}\,\binom{2\ell}{\ell-k}\,
    U_{2k-1}\left(\frac{s}{2R}\right)\Bigg]
\end{split}
\end{equation}
where
$U_{k}(x)=\sin\big((k+1)\arccos(x)\big)/\sin\big(\arccos(x)\big)$ are
the Chebyshev polynomials of the second kind \cite[Ch.~22]{AS}. For
odd dimension, $m=2\ell+1$, one obtains the following expression
\begin{equation}
   f^{}_{2\ell+1}(s) \; = \; 
    1 - 
    \frac{\Gamma(\ell+\tfrac{3}{2})}{2^{2\ell-1}\sqrt{\pi}\,\Gamma(\ell+1)}\,
    \sum\limits_{k=0}^{\ell}\frac{(-1)^{k}}{2k+1}\,\binom{2\ell+1}{\ell-k}\,
    T_{2k+1}\left(\frac{s}{2R}\right)
\label{eq:fodd}
\end{equation}
in terms of the Chebyshev polynomials
$T_{k}(x)=\cos\big(k\arccos(x)\big)$ of the first kind
\cite[Ch.~22]{AS}.  Eqs.~(\ref{eq:fgen})--(\ref{eq:fodd}) are valid
for $0\le s<2R$; for distances larger than the diameter, the overlap
vanishes, hence $f^{}_{m}(s)=0$ for $s\ge 2R$. Eq.~(\ref{eq:f2}) is
recovered from (\ref{eq:feven}) for $\ell=1$. The functions
$f^{}_{m}(s)$ for various dimensions $m$ are shown in
Fig.~\ref{fig:rad}. Unfortunately, for $m>2$, there is no simple
generalisation of Eq.~(\ref{a-c}), because the $*$-map is then more
complicated.\smallskip

\begin{figure}
\centerline{\epsfysize=0.45\textheight\epsfbox{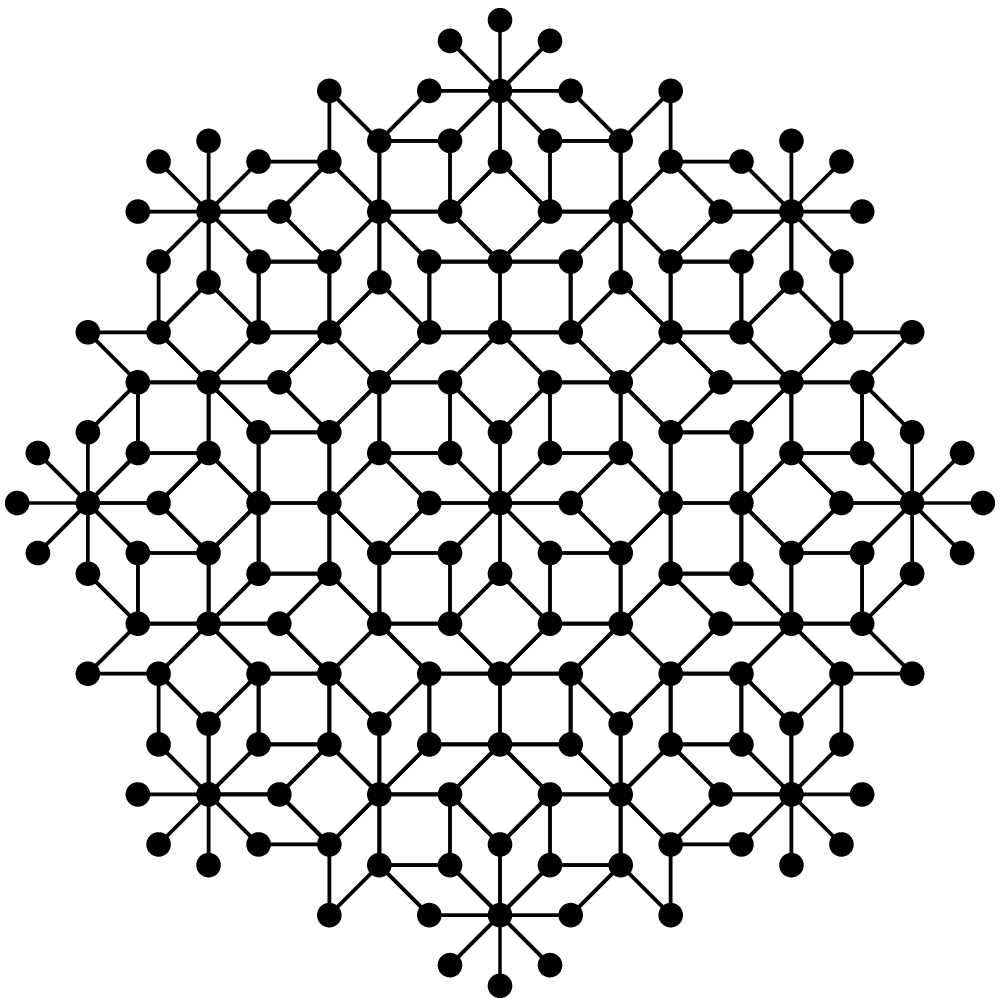}\hspace{0.05\textwidth}
\epsfysize=0.45\textheight\epsfbox{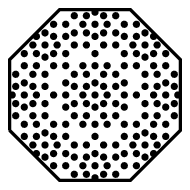}}
\caption{A patch of the Ammann-Beenker tiling with vertex set
$\varLambda_{\rm AB}$ (left) and the $*$-image of $\varLambda_{\rm AB}$ inside
the octagonal window in internal space (right), with relative scale as
described in the text.
\label{fig:ab}}
\end{figure}

Let us finally consider an eightfold symmetric model set in the plane,
based on the classical Ammann-Beenker or octagonal tiling, compare
\cite{AGS,BGJR,BJ} and references therein. It is usually described by
projection from four dimensions, where we use the lattice
$\varGamma=\sqrt{2}\,\ZZ^4$. The projections $\pi$ and $\pi_{H}^{}$ of
(\ref{candp}) are essentially determined by compatibility with
eightfold symmetry. In a convenient coordinatisation \cite{BJ}, the
images $\bs{a}_{j}$, $j\in\{1,2,3,4\}$, of the standard basis vectors
of the lattice have unit length in physical space, and the same is
true of the corresponding projections $\bs{a}_{j}^{*}$ in internal
space. Observing that $\ZZ\bs{a}_{1}+\ZZ\bs{a}_{2}+
\ZZ\bs{a}_{3}+\ZZ\bs{a}_{4} = \ZZ[\xi_{8}^{}]$, we can continue with a
formulation based on the cyclotomic integers, compare \cite{P}.  For
the Ammann-Beenker tiling, the window is then a regular octagon $O$ of
unit edge length, see Fig.~\ref{fig:ab}. Note that the window is
invariant under the symmetry group $D_8$ of order $16$.

\begin{figure}[t]
\centerline{\epsfxsize=0.7\textwidth\epsfbox{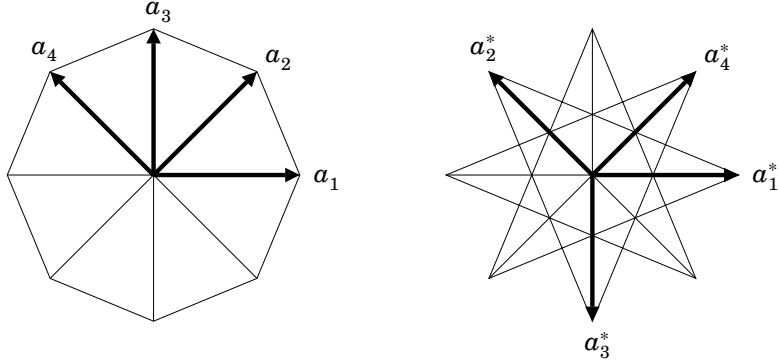}}
\caption{Vectors $\bs{a}_{j}$ in physical and $\bs{a}_{j}^{*}$
in internal space related by the $*$-map.\label{fig:star}}
\end{figure}

Explicitly, the corresponding point set in the plane is given by
\begin{equation}
\varLambda_{\text{AB}}\; =\;
\big\{z\in\ZZ[\xi_{8}]\mid z^{*}\in O\big\},
\end{equation}
where ${\,}^*$ is the Galois automorphism defined by
$\xi_{8}^{}\mapsto\xi_{8}^{3}$. If we choose $\xi_{8}^{}=\xi=e^{2\pi
i/8}$ and identify $\RR^{2}$ with $\CC$, this gives
$\bs{a}_{j}=\xi^{j-1}$, $1\le j\le 4$, while the $*$-images satisfy
$\bs{a}_{j}^{*}=\xi^{3(j-1)}$, compare Fig.~\ref{fig:star}.

\begin{figure}[b]
\centerline{\epsfxsize=0.6\textwidth\epsfbox{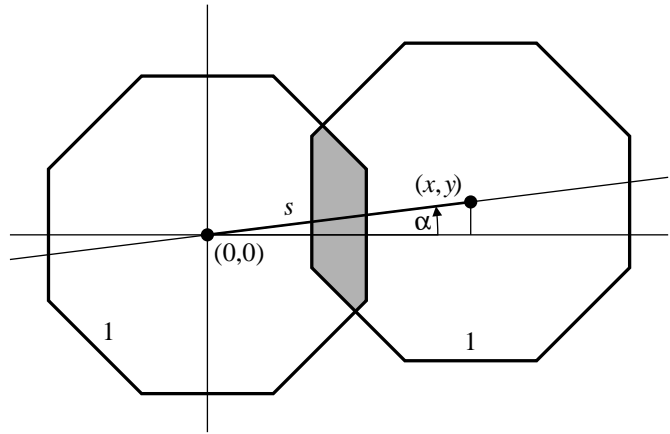}}
\caption{Two overlapping regular octagons at distance
$s$.\label{fig:oct}}
\end{figure}

\begin{figure}
\centerline{\epsfxsize=0.52\textwidth\epsfbox{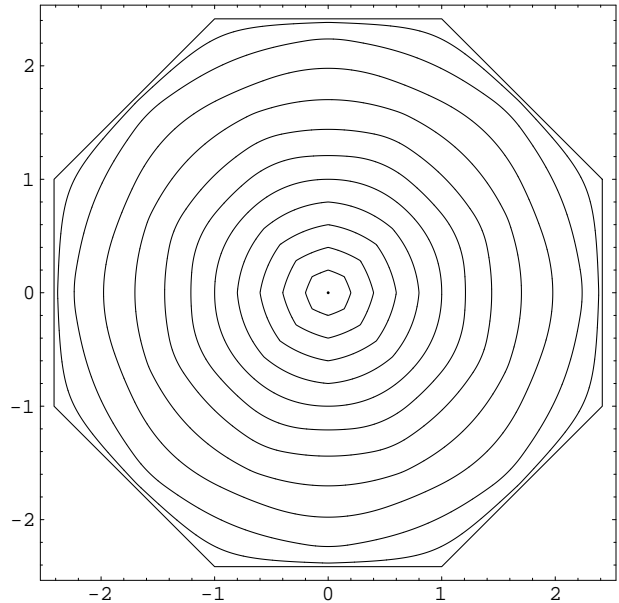}}
\caption{Contour map of the function $f^{}_{\text{AB}}(s,\alpha)$ of
(\ref{fab}) in internal space. The contours show the decrease of
$f^{}_{\text{AB}}$ from its maximum value $f^{}_{\text{AB}}(0,0)=1$ to
$f^{}_{\text{AB}}(s,\alpha)=0$ outside the outer octagonal
contour.\label{fig:fab}}
\end{figure}

A somewhat tedious, but elementary calculation on the basis
of  Fig.~\ref{fig:oct} gives
\begin{fact}\label{covar}
The covariogram of the regular octagon of edge length one is
\begin{equation}
f^{}_{\rm AB}(s,\alpha)\; =\; 
\begin{cases}
0 & \mbox{if\/ $\lambda\le x$}\\
\frac{(\lambda-2)(x+y)x}{2}+\frac{\lambda(1-x)}{2}+\frac{(1-y)}{2}\quad
& \mbox{if\/ $\lambda-y\le x\le\lambda$}\\
\frac{(\lambda-2)(x^2-y^2)}{4}-\frac{(\lambda-1)x}{2}+\frac{\lambda+2}{4}\quad
& \mbox{if\/ $1+y\le x\le\lambda-y$}\\
\frac{(\lambda-2)(x-y-1)y}{2}-\frac{x}{2}+1
& \mbox{if\/ $0\le x\le 1+y$}
\end{cases}\label{fab}
\end{equation}
where\/ $\lambda = 1+\sqrt{2}$,
\[
x \; = \; s\,\cos\alpha'  \; \in \;
        \Big[{\textstyle\frac{s}{2}\sqrt{2+\sqrt{2}}},s\Big]\, ,\qquad
y \; = \; s\,\sin\alpha'  \; \in \; 
\Big[0,{\textstyle\frac{s}{2}\sqrt{2-\sqrt{2}}}\,\Big] \, ,
\]
and where\/ $\alpha'$ is the unique angle in the interval\/
$[0,\frac{\pi}{8}]$ that is related to\/ $\alpha$ by the\/ $D_8$
symmetry of the octagon.\hfill\qed
\end{fact}
A contour map of $f^{}_{\text{AB}}(s,\alpha)$ is shown in
Fig.~\ref{fig:fab}. It demonstrates that the previous consideration of
a circular window is actually a reasonable approximation to this case.
It is sufficient for most applications concerning (powder)
diffraction, compare \cite[Ch.~3]{JS}.

We can now calculate the averaged shelling coefficient $a(r^2)$ of
(\ref{ashell1}) explicitly for any distance $r$ in $\varLambda_{\rm
AB}$.  The results for all distances with $0<r^{2}\le 5$ are
summarised in Table~\ref{abtab}. They confirm the results of
\cite{BGJR} which had been obtained numerically.
 
As an explicit example, let us consider the shortest distance in the
model set.  This is $r=\sqrt{2-\sqrt{2}}=2\sin(\frac{\pi}{8})$ which
is realised by the short diagonal of the rhomb. In this case, there
are eight numbers $z\in\varLambda_{\text{AB}}-\varLambda_{\text{AB}}$
that contribute to Eq.~(\ref{ashell1}). They form a single $D_8$
orbit; one representative is listed in Table~\ref{abtab}.  Due to the
symmetry of the window, the contribution of each member of the orbit
is the same, so it suffices to consider a representative and to
multiply the result by the corresponding orbit length. Choosing
$z=1-\xi$, we find $z^{*}=1-\xi^{3}$, hence
$|z^{*}|=\sqrt{2+\sqrt{2}}=2\cos(\frac{\pi}{8})$; the corresponding
angle is $\alpha=-\frac{\pi}{8}$, hence $\alpha'=\frac{\pi}{8}$. The
distance in internal space is rather large and the overlap area
correspondingly small. Fact~\ref{covar} yields
$x=1+y=\big(2+\sqrt{2}\,\big)/2$; hence the coefficient is given by
\begin{equation}
   a\big(2-\sqrt{2}\,\big)\; =\; 
   8\, f^{}_{\rm AB}\big(2\cos{(\tfrac{\pi}{8})},\tfrac{\pi}{8}\big)
   \;=\; 8\big(1-\tfrac{x}{2}\big)
   \;=\; 4-2\sqrt{2},
\end{equation}
compare Table~\ref{abtab}.

The other entries in Table~\ref{abtab} are calculated along the same
lines. Note that $s$ can be calculated from $r$ directly via
$s^2={(r^{2})}^*$, where ${}^*$ coincides with algebraic conjugation
in $\QQ(\sqrt{2}\,)$, defined by $\sqrt{2}\mapsto -\sqrt{2}$.
Continuing the calculations, one faces increasing complication with
growing distance $r$, and, in general, one has to expect contributions
from several $D_8$ orbits. For $r=\sqrt{3}$, there is still only a
single orbit, this time of length $16$. Hence it again suffices to
consider a single representative, for instance $z=1-\xi-\xi^{2}$ whose
$*$-image is
$z^{*}=\tfrac{1}{2}\big((2+\sqrt{2}\,)+i\,(2-\sqrt{2}\,)\big)$.  More
generally, the standard orbit analysis reduces the sum in
(\ref{ashell1}) to a formula with one contribution per $D_8$-orbit,
weighted with the corresponding orbit length. The latter is $8$
whenever the corresponding angle $\alpha$ of Fact~\ref{covar} is an
integer multiple of $\tfrac{\pi}{8}$ (which corresponds to symmetry
directions of the octagon), and $16$ otherwise.

\begin{table}
\caption{The averaged shelling numbers for distances $r$ with
$0<r^2\le 5$ in $\varLambda_{\rm AB}$.
Representatives $z$ are given in terms of $\xi=e^{2\pi i/8}$.  The
examples listed comprise a single $D_8$ orbit each.
\label{abtab}}
\renewcommand{\arraystretch}{2}
\begin{tabular}{cccccc}
\hline
$r$ &
$z$ &
orbit length &
$s$ &
$\alpha'$ &
$a(r^2)$\\
\hline
$\sqrt{2-\sqrt{2}}$ &
$1-\xi$ &
$8$ &
$\sqrt{2+\sqrt{2}}$ & $\frac{\pi}{8}$ & $4-2\sqrt{2}$\\
$1$ & 
$1$ &
$8$ &
$1$ & $0$ & $4$\\
$\sqrt{2}$ &
$1+\xi^2$ &
$8$ & 
$\sqrt{2}$ & $0$ & $6\sqrt{2}-6$\\
$\sqrt{3}$ &
$1-\xi-\xi^{2}$ &
$16$ &
$\sqrt{3}$ & $\arctan\big(\frac{2-\sqrt{2}}{2+\sqrt{2}}\big)$ 
& $20-12\sqrt{2}$ \\
$\sqrt{2+\sqrt{2}}$ &
$1+\xi$ &
$8$ &
$\sqrt{2-\sqrt{2}}$ & $\frac{\pi}{8}$ & $36-22\sqrt{2}$\\
$2$ &
$2$ &
$8$ &
$2$ & $0$ & $2\sqrt{2}-2$\\ 
$\sqrt{5}$ &
$2+\xi^2$ &
$16$ & 
$\sqrt{5}$ &
$\arctan\big(\frac{1}{3}\big)$ &
$40\sqrt{2}-56$\\
\hline
\end{tabular}
\renewcommand{\arraystretch}{1}
\end{table}
 
The averaged shelling numbers $a(r^2)$ of Table~\ref{abtab}, as well
as the numerically determined values in \cite[Tab.~2]{BGJR}, are
always elements of $2\hspace{0.5pt}\ZZ[\sqrt{2}\,]$, i.e., numbers of
the form $2k+2\ell\sqrt{2}$ with $k,\ell\in\ZZ$. The formula
(\ref{ashell2}) for the covariogram of the regular octagon only
implies $a(r^2)$ to be in $\QQ(\sqrt{2}\,)$. However,
Eqs.~(\ref{ashell1}) and (\ref{weyl}) show that the averaged shelling
number $a(r^2)$ is expressed as a finite sum of cluster
frequencies. The latter are elements of the so-called {\em frequency
module}\/ $F_{\rm AB}$ of the Ammann-Beenker model set
$\varLambda_{\rm AB}$, i.e., the integer span of the frequencies of
finite clusters of $\varLambda_{\rm AB}$ (i.e., of all intersections
$\varLambda_{\rm AB}\cap C$ with $C\subset\RR^2$ compact). Since
$\varLambda_{\rm AB}$ is a Delone set of finite local complexity
\cite{Lag,P}, there are only countably many different clusters up to
translations, and only finitely many up to a given diameter.

The frequency module has originally been calculated by means of
$C^{*}$-algebraic $K$-theory \cite{Belli,FHK}, but can also be
obtained from simpler cohomological considerations \cite{Franz}. For
our case, the result can be simplified, by a short calculation, as
follows.
\begin{fact}\cite{Belli,Franz}
The frequency module of the Ammann-Beenker vertex set $\varLambda_{\rm
AB}$ is 
\[
F_{\rm AB}\;=\;\Big\{ \frac{k+\ell\sqrt{2}}{8}\;\big|\; k,\ell\in\ZZ,\;
k+\ell \mbox{ even}\Big\}.
\]
Thus, $8F_{\rm AB}$ is an index\/ $2$ submodule of\/ $\ZZ[\sqrt{2}\,]$.
\hfill\qed
\end{fact}
Since the window of $\varLambda_{\rm AB}$ is a regular octagon, each
cluster that contributes to the sum in Eq.~(\ref{ashell1}) occurs in
either $8$ or $16$ orientations with the same frequency. Consequently,
the averaged shelling numbers $a(r^{2})$ are elements of $8F_{\rm
AB}$, though we presently do not know whether they generate the full
frequency module or a subset thereof. This is consistent with the
findings of Table~\ref{abtab} and \cite[Tab.~2]{BGJR}, and establishes
an interesting link between geometric and topological properties of
model sets \cite{AP,FHK}.
\medskip

A similar treatment is possible for other examples, such as the
tenfold symmetric rhombic Penrose \cite{BGJR} and T\"{u}bingen
triangle \cite{MU1} tilings, or the twelvefold symmetric shield tiling
\cite{MU2}.  However, it is clear that, for the standard model sets,
the calculation becomes rather involved, and we presently do not know
how to determine the averaged shelling function via a generating
function such as (\ref{theta2}) in closed form, or even whether that
is the most promising way to proceed. For physical applications,
however, often the first few terms are sufficient, and they can be
calculated exactly from the projection method, see \cite{BGJR,MU1,MU2}
for some tables and further results.

\section*{Acknowledgments}

It is a pleasure to thank Robert V.\ Moody and Alfred Weiss for
helpful discussions and suggestions. This work was partially supported
by the German Research Council (DFG). We also thank two anonymous
referees for constructive and helpful comments, and the Erwin
Schr\"{o}dinger International Institute for Mathematical Physics in
Vienna for support during a stay in winter 2002/2003 where the revised
version was prepared.


\bigskip

\begin{thebibliography}{99}
\small

\bibitem{AS}
M.~Abramowitz and I.~A.~Stegun (eds.),
{\em Pocketbook of Mathematical Functions},
Harri Deutsch, Thun (1984).

\bibitem{AGS}
R.~Ammann, B.~Gr\"{u}nbaum and G.~C.~Shephard,
Aperiodic tiles,
{\em Discr.\ Comput.\ Geom.}\/ {\bf 8} (1992) 1--25.

\bibitem{AP}
J.~E.~Anderson and I.~F.~Putnam,
Topological invariants for substitution tilings and their associated 
$C^*$-algebras,
{\em Ergodic Theory Dynam.\ Systems}\/ {\bf 18} (1998) 509--537.

\bibitem{MB}
M.~Baake,
A guide to mathematical quasicrystals,
in: {\em Quasicrystals}, eds.\ J.-B.\ Suck,  M.\ Schreiber and P.\
H\"{a}ussler, Springer, Berlin (2002), pp.~17--48; 
math-ph/9901014.

\bibitem{BGlit}
M.~Baake and U.~Grimm,
A guide to quasicrystal literature,
in: \cite{BMbook}, pp.~371--373.

\bibitem{MU1}
M.~Baake and U.~Grimm,
Quasicrystalline combinatorics,
to appear in: {\em Proceedings of GROUP\/$24$},
eds.\ J.\ P.\ Gazeau and R.\ Kerner, 
IOP, Bristol (2003);
preprint mp\_arc/02-392.

\bibitem{MU2}
M.~Baake and U.~Grimm,
Combinatorial problems of (quasi-)crystallography,
to appear in: {\em Quasicrystals -- Structure and Physical Properties}, 
ed.\ H.-R. Trebin, Wiley-VCH, Berlin (2003); 
preprint math-ph/0212015.

\bibitem{BGJR}
M.~Baake, U.~Grimm, D.~Joseph and P.~Repetowicz,
Averaged shelling for quasicrystals,
{\em Mat.\ Sci.\ Eng.\ A}\/ {\bf 294--296} (2000) 441--445;
math.MG/9907156.

\bibitem{BJ}
M.~Baake and D.~Joseph,
Ideal and defective vertex configurations in the planar octagonal 
quasilattice,
{\em Phys.\ Rev.\ B}\/ {\bf 42} (1990) 8091--8102.

\bibitem{BJS}
M.~Baake, D.~Joseph and M.~Schlottmann,
The root lattice $D_4$ and planar quasilattices with octagonal and
dodecagonal symmetry,
{\em Int.\ J.\ Mod.\ Phys.\ B}\/ {\bf 5} (1991) 1927--1953.

\bibitem{BMbook}
M.~Baake and R.~V.~Moody (eds.),
{\em Directions in Mathematical Quasicrystals},
CRM Monograph Series, vol.~13,
AMS, Providence, RI (2000).

\bibitem{BM2001}
M.~Baake and R.~V.~Moody,
Weighted Dirac combs with pure point diffraction,
preprint math.MG/0203030.

\bibitem{Belli}
J.~Bellissard, E.~Contensou and A.~Legrand,
$K$-th\'{e}orie des quasi-cristaux, image par la trace: 
Le cas du r\'{e}seau octogonal,
{\em C.\ R.\ Acad.\ Sci.\ Paris S\'{e}r.\ I Math.}\/ 
{\bf 326} (1998) 197--200.

\bibitem{BS}
K.~B\"{o}r\"{o}czky Jr.\ and U.~Schnell,
Quasicrystals and the Wulff-shape,
{\em Discr.\ Comput.\ Geom.}\/ {\bf 21} (1999) 421--436.

\bibitem{BSW}
K.~B\"{o}r\"{o}czky Jr., U.~Schnell and J.~M.~Wills,
Quasicrystals, parametric density, and Wulff-shape,
in: \cite{BMbook}, pp.~259--276.

\bibitem{CS}
J.~H.~Conway and N.~J.~A.~Sloane,
{\em Sphere Packings, Lattices and Groups}, 3rd ed.,
Springer, New York (1999).

\bibitem{DFKPRW}
M.~Daberkow, C.~Fieker, J.~Kl\"{u}ners, M.~Pohst, K.~Roegner and
K.~Wildanger,
KANT V4,
{\em J.\ Symbolic Comp.}\/ {\bf 24} (1997) 267--283.

\bibitem{ES}
V.~Elser and N.~J.~A.~Sloane,
A highly symmetric four-dimensional quasicrystal,
{\em J.\ Phys.\ A: Math.\ Gen.}\/ {\bf 20} (1987) 6161--6168.

\bibitem{FHK}
A.~H.~Forrest, J.~R.~Hunton and J.~Kellendonk, 
{\em Topological Invariants for Projection Method Patterns}, 
Memoirs of the AMS, vol.~159,
AMS, Providence, RI (2002).

\bibitem{Franz}
F.~G\"{a}hler,
private communication (2002).

\bibitem{Gar}
D.~A.~Garbanati,
Units with norm $-1$ and signatures of units,
{\em J.\ Mathematik $($Crelle\/$)$}\/ {\bf 283/284} (1976) 164--175.

\bibitem{Gneit}
T.~Gneiting,
Radial positive definite functions generated by Euclid's hat,
{\em J.\ Multivariate Anal.}\/ {\bf 69} (1999) 88--119.

\bibitem{Hof}
A.~Hof,
Uniform distribution and the projection method,
in: \cite{Pat}, pp.~201--206.

\bibitem{Hux}
M.~N.~Huxley,
{\em Area, Lattice Points, and Exponential Sums},
Clarendon Press, Oxford (1996).

\bibitem{JS}
R.~Jenkins and R.~L.~Snyder,
{\em Introduction to X-ray Powder Diffractometry},
Wiley, New York (1996).

\bibitem{kant}
KANT/KASH,
{\verb+http://www.math.tu-berlin.de/~kant/kash.html+}.

\bibitem{Lag}
J.~C.~Lagarias,
Geometric models for quasicrystals I. Delone sets of finite type,
{\em  Discr.\ Comput.\ Geom.}\/ {\bf 21} (1999) 161--191.

\bibitem{LM}
J.-Y.~Lee and R.~V.~Moody,
Lattice substitution systems and model sets,
{\em Discr.\ Comput.\ Geom.}\/ {\bf 25} (2001) 173--201;
math.MG/0002019.

\bibitem{Mermin}
N.~D.~Mermin, D.~S.~Rokhsar and D.~C.~Wright,
Beware of $46$-fold symmetry: The classification of two-dimensional
quasicrystallographic lattices,
{\em Phys.\ Rev.\ Lett.}\/ {\bf 58} (1987) 2099--2101.

\bibitem{Mbook}
R.~V.~Moody (ed.),
{\em The Mathematics of Long-Range Aperiodic Order},
NATO ASI Series C 489,
Kluwer, Dordrecht (1997).

\bibitem{Moody}
R.~V.~Moody,
Model sets: A survey, in: 
{\em From Quasicrystals to More Complex Systems}, 
eds.\ F.\ Axel, F.\ D\'enoyer and J.\ P.\ Gazeau,
EDP Sciences, Les Ulis, and
Springer, Berlin (2000), pp.\ 145--166;
math.MG/0002020.

\bibitem{Moody2001}
R.~V.~Moody,
Uniform distribution in model sets,
{\em Can.\ Math.\ Bulletin}\/ {\bf 45} (2002)
123--130.

\bibitem{MW}
R.~V.~Moody and A.~Weiss,
On shelling $E_8$ quasicrystals,
{\em J.\ Number Th.}\/ {\bf 47} (1994) 405--412.

\bibitem{MS}
J.~Morita and K.~Sakamoto,
Octagonal quasicrystals and a formula for shelling,
{\em J.\ Phys.\ A: Math.\ Gen.}\/ {\bf 31} (1998) 9321--9325.

\bibitem{MS2}
J.~Morita and K.~Sakamoto,
Shell structure of dodecagonal quasicrystals associated with root 
system $F_4$ and cyclotomic field $Q(\zeta(12))$,
{\em Commun.\ Algebra}\/ {\bf 28} (2000) 255--263.

\bibitem{Neu}
J.~Neukirch,
{\em Algebraic Number Theory},
Springer, Berlin (1999).

\bibitem{Pat}
J.~Patera (ed.),
{\em Quasicrystals and Discrete Geometry},
Fields Institute Monographs, vol.\ 10,
AMS, Providence, RI (1998).

\bibitem{P}
P.~A.~B.~Pleasants,
Designer quasicrystals: Cut-and-project sets with pre-assigned 
properties,
in: \cite{BMbook}, pp.~95--141.

\bibitem{PBR}
P.~A.~B.~Pleasants, M.~Baake and J.~Roth,
Planar coincidences for $N$-fold symmetry,
{\em J.\ Math.\ Phys.}\/ {\bf 37} (1996) 1029--1058.

\bibitem{SM}
J.-F.~Sadoc and R.~Mosseri,
The $E_8$ lattice and quasicrystals: Geometry, number theory,
and quasicrystals, 
{\em J.\ Phys.\ A: Math.\ Gen.}\/ {\bf 26} (1993) 1789--1809.

\bibitem{Martin1}
M.~Schlottmann,
Cut-and-project sets in locally compact Abelian groups,
in: \cite{Pat}, pp.~247--264.

\bibitem{Martin2}
M.~Schlottmann,
Generalized model sets and dynamical systems,
in: \cite{BMbook}, pp.~143--159.

\bibitem{S}
M.~Senechal,
{\em Quasicrystals and Geometry},
Cambridge University Press, Cambridge (1995); 
corrected reprint (1996).

\bibitem{Solo}
B.~Solomyak,
Spectrum of dynamical systems arising from Delone sets,
in \cite{Pat}, pp.~265--275.

\bibitem{Wash}
L.~C.~Washington,
{\em Introduction to Cyclotomic Fields}, 2nd ed.,
Springer, New York (1997).

\bibitem{Al}
A.~Weiss,
On shelling icosahedral quasicrystals,
in: \cite{BMbook}, pp.~161--176.

\end{thebibliography}
\end{document}